\documentclass[MSNbibl,nameyear,seceqn,dvips]{arxstspdf}
\usepackage{flushend}
\usepackage{stfloats}
\usepackage{graphicx}
\volume{27}
\issue{4}
\pubyear{2012}
\firstpage{481}
\lastpage{499}
\doi{10.1214/12-STS392} 

\makeatletter

\newcommand{\bb}{\mathbf{b}}
\newcommand{\rP}{\mathrm{P}}
\newcommand{\bfr}{\mathbf{r}}
\newcommand{\bu}{\mathbf{u}}
\newcommand{\bv}{\mathbf{v}}
\newcommand{\bx}{\mathbf{x}}
\newcommand{\tX}{\widetilde{X}}
\newcommand{\by}{\mathbf{y}}
\newcommand{\tby}{\tilde{\mathbf{y}}}
\newcommand{\bz}{\mathbf{z}}
\newcommand{\tbz}{\tilde{\mathbf{z}}}
\newcommand{\bzero}{\mathbf{0}}
\newcommand{\veps}{\varepsilon}
\newcommand{\argmin}{\operatorname{argmin}}
\newcommand{\real}{\mathbb{R}}

\newcommand{\htheta}{\hat{\theta}}
\newcommand{\btheta}{\bolds{\theta}}
\newcommand{\bbeta}{\bolds{\beta}}
\newcommand{\hbbeta}{\hat{\bolds{\beta}}}
\newcommand{\tbbeta}{\tilde{\bolds{\beta}}}
\newcommand{\gam}{\gamma}
\newcommand{\lam}{\lambda}
\newcommand{\drho}{\dot{\rho}}
\newcommand{\hbtheta}{\widehat{\bolds{\theta}}}
\newcommand{\bveps}{\bolds{\varepsilon}}

\newproclaim{definition}{Definition}[section]
\newtheorem{theorem}{Theorem}[section]
\newtheorem{lemma}{Lemma}[section]
\newtheorem{corollary}{Corollary}[section]

\newproclaim{example}{Example}[section]
\newproclaim{remark}{Remark}[section]
\makeatother

\begin{document}
\begin{frontmatter}

\title{A Selective Review of Group Selection in High-Dimensional Models}
\runtitle{Group Selection}

\begin{aug}
\author[a]{\fnms{Jian} \snm{Huang}\corref{}\ead[label=e1]{jian-huang@uiowa.edu}},
\author[b]{\fnms{Patrick} \snm{Breheny}\ead[label=e2]{patrick.breheny@uky.edu}}
\and
\author[c]{\fnms{Shuangge} \snm{Ma}\ead[label=e3]{shuangge.ma@yale.edu}}
\runauthor{J. Huang, P. Breheny and S. Ma}

\affiliation{University of Iowa, University of Kentucky and Yale University}

\address[a]{Jian Huang is Professor, Department of Statistics and
Actuarial Science, 241 SH,
University of Iowa, Iowa City, Iowa 52242, USA  \printead{e1}.}
\address[b]{Patrick Breheny is Assistant Professor, Department of Statistics,
University of Kentucky, Lexington, Kentucky 40506, USA \printead{e2}.}
\address[c]{Shuangge Ma is Associate Professor, Division of
Biostatistics, School of
Public Health, Yale University, New Haven, Connecticut 06520, USA
\printead{e3}.}

\end{aug}

%
\begin{abstract}
Grouping structures arise naturally in many statistical
modeling problems. Several methods have been proposed for variable
selection that respect
grouping structure in variables. Examples include the group LASSO and
several concave group selection methods.
In this article, we give a selective review of group selection
concerning methodological developments,
theoretical properties and computational algorithms. We pay particular
attention to group selection
methods involving concave penalties. We address both group selection
and bi-level selection methods.
We describe several applications of these methods in nonparametric
additive models, semiparametric
regression, seemingly unrelated regressions, genomic data analysis and
genome wide association studies.
We also highlight some issues that require further study.
\end{abstract}

%

\end{frontmatter}
%

\section{Introduction}
Consider a linear regression model with $p$ predictors. Suppose the
predictors can be naturally divided into $J$ nonoverlapping groups, and
the model is written as
%
\begin{equation}\label{Mod1}
\by=\sum_{j=1}^{J}X_{j}\bbeta_{j}+\bveps,
\end{equation}
where $\by$ is an $n
\times1$ vector of response variables, $X_{j}$ is the $n \times d_j$
design\vspace*{1pt} matrix of the $d_j$ predictors in the $j$th group,
$\bbeta_{j}=(\beta_{j1}, \ldots, \beta_{jd_j})'\in\real^{d_j}$ is the
$d_j \times1$ vector of regression coefficients of the $j$th group and
$\bveps$ is the error vector. Without loss of generality, we take both
the predictors and response to be centered around the mean. It is
desirable to treat each group of variables as a unit and take advantage
of the grouping structure present in these models when estimating
regression coefficients and selecting important variables.

Many authors have considered the problem of group selection in various
statistical modeling problems. Bakin (\citeyear{Bakin}) proposed the group LASSO
and a computational algorithm. This method and related group selection
methods and algorithms were further developed by Yuan and Lin (\citeyear{YuanL2006}).
The group LASSO uses an $\ell_2$ norm of the coefficients associated
with a group of variables in the penalty function and is a natural
extension of the LASSO (Tibshirani, \citeyear{Tibshirani96}). Antoniadis and Fan (\citeyear{AnFan})
studied a class of block-wise shrinkage approaches for regularized
wavelet estimation in nonparametric regression problems. They discussed
several ways to shrink wavelet coefficients in their natural blocks,
which include the blockwise hard- and soft-threshold rules. Meier, van
de Geer and B\"{u}hlmann (\citeyear{MGB2008}) studied the group LASSO for logistic
regression. Zhao, Rocha and Yu (\citeyear{ZhaoRY2009}) proposed a quite general
composite absolute penalty for group selection, which includes the
group LASSO as a special case. Huang, Ma, Xie and Zhang (\citeyear{HuangMXZ})
considered the problem of simultaneous group and individual variable
selection, or bi-level selection, and proposed a group bridge method.
Breheny and Huang (\citeyear{BrehenyHuang1}) proposed a general framework for bi-level
selection in generalized linear models and derived a local coordinate
descent algorithm.

Grouping structures can arise for many reasons, and give rise to quite
different modeling goals. Common examples include the representation
of multilevel categorical covariates in a regression model by a group
of indicator variables, and the representation of the effect of a
continuous variable by a set of basis functions. Grouping can also be
introduced into a model in the hopes of taking advantage of prior
knowledge that is scientifically meaningful. For example, in gene
expression analysis, genes belonging to the same biological pathway can
be considered a group. In genetic association studies, genetic markers
from the same gene can be considered a group. It is desirable to take
into account the grouping structure in the analysis of such data.

Depending on the situation, the individual variables in the groups may
or may not be meaningful scientifically. If they are not, we are
typically not interested in selecting individual variables; our
interest is entirely in {\em group selection}. However, if individual
variables are meaningful, then we are usually interested in selecting
important variables as well as important groups; we refer to this as
{\em bi-level selection}. For example, if we represent a continuous
factor by a set of basis functions, the individual variables are an
artificial construct, and selecting the important members of the group
is typically not of interest. In the gene expression and genetic
marker examples, however, selection of individual genes/markers is just
as important as selecting important groups. In other examples, such as
a group of indicator functions for a categorical variable, whether we
are interested in selecting individual members depends on the context
of the study.

We address both group selection and bi-level selection in this review.
The distinction between these two goals is crucial for several reasons.
Not only are different statistical methods used for each type of
problem, but as we will see, the predictors in a group can be made
orthonormal in settings where bi-level selection is not a concern.
This has a number of ramifications for deriving theoretical results and
developing algorithms to fit these models.

We give a selective review of group selection concerning methodological
developments, theoretical\break properties and computational algorithms. We
describe several important applications of group selection
and bi-level selection in nonparametric additive models, semiparametric
regression, seemingly unrelated regressions, genomic data analysis and
genome wide association studies. We also highlight some issues that
require further study. For the purposes of simplicity, we focus
on penalized versions of least squares regression in this review. Many
authors have extended these models to other loss functions, in
particular those of generalized linear models. We attempt to point out
these efforts when relevant.

\section{Group Selection Methods}\label{sec2}

\subsection{Group LASSO}
For a column vector $\bv\in\real^d$ with $d \ge1$ and a positive
definite matrix $R$, denote $\|\bv\|_2=(\bv'\bv)^{1/2}$ and
$\|\bv\|_{R} = (\bv'R \bv)^{1/2}$. Let $\bbeta=(\bbeta_1', \ldots,
\bbeta_J')'$, where $\bbeta_j \in\real^{d_j}$. The group LASSO
solution $\hbbeta(\lam)$ is defined as a minimizer of
%
\begin{equation}\label{gL1}
\frac{1}{2n}\Biggl\|\by-\sum_{j=1}^J X_j \bbeta_j\Biggr\|_2^2 + \lam\sum_{j=1}^J
c_j \|\bbeta_j\|_{R_j},
\end{equation}
where $\lam\ge0$ is the penalty
parameter and $R_j$'s are $d_j \times d_j$ positive definite matrices.
Here the $c_j$'s in the penalty are used to adjust for the group sizes.
A reasonable choice is $c_j = \sqrt{d_j}$. Because (\ref{gL1}) is
convex, any local minimizer of (\ref{gL1}) is also a global minimizer
and is characterized by the Karush--Kuhn--Tucker conditions as given in
Yuan and Lin (\citeyear{YuanL2006}). It is possible, however, for multiple solutions
to exist, as (\ref{gL1}) may not be strictly convex in situations where
the ordinary least squares estimator is not uniquely defined.

An important question in the definition of group LASSO is the choice of
$R_j$. For orthonormal $X_j$ with $X_j'X_j/n=I_{d_j}$, $j=1, \ldots,
J$, Yuan and Lin (\citeyear{YuanL2006}) suggested taking $R_j= I_{d_j}$. However,
using $R_j=I_{d_j}$ may not be appropriate, since the scales of the
predictors may not be the same. In general, a reasonable choice of
$R_j$ is to take the Gram matrix based on $X_j$, that is, $R_j =
X_j'X_j/n$, so that the penalty is proportional to $\|X_j\bbeta_j\|_2$.
This is equivalent to performing standardization at the group level,
which can be seen as follows. Write $R_j = U_j'U_j$ for a $d_j\times
d_j$ upper triangular matrix $U_j$ via Cholesky decomposition. Let
$\tX_j=X_jU_j^{-1}$ and $\bb_j = U_j\bbeta_j$. Criterion (\ref{gL1})
becomes
%
\begin{equation}\label{gL2} \frac{1}{2n}\Biggl\|\by-\sum_{j=1}^J \tX_j
\bb_j\Biggr\|
_2^2+\lam
\sum_{j=1}^{J} c_j \|\bb_j\|_2.
\end{equation}
The solution to the original
problem (\ref{gL1}) can be obtained by using the transformation
$\bbeta_j =U_j^{-1} \bb_j$. By the definition of $U_j$, we have $n^{-1}\tX_j'\tX_j = I_{d_j}$.
Therefore, by using this choice of
$R_j$, without loss of generality, we can assume that $X_j$ satisfies
$n^{-1}X_j'X_j=I_{d_j}, 1\le j \le J$. Note that we do not assume
$X_j$ and $X_k$, $j\neq k$, are orthogonal.

The above choice of $R_j$ is easily justified in the special case where
$d_j=1, 1\le j \le J$. In this case, the group LASSO simplifies to the
standard LASSO and $R_j=\|X_j\|^2/n$ is proportional to the sample
variance of the $j$th predictor. Thus, taking $R_j$ to be the Gram
matrix is the same as standardizing the predictors before the analysis,
which is often recommended when applying LASSO for variable selection.

Several authors have studied the theoretical properties of the group
LASSO, building on the ideas and approaches for studying the behavior
of the LASSO, on which there is an extensive literature; see
B\"{u}hlmann and van de Geer (\citeyear{BG2011}) and the references therein. Bach
(\citeyear{Bach2008}) showed that the group LASSO is group selection consistent in a
random design model for fixed $p$ under a variant of the
irrepresentable condition (Meinshausen and B\"{u}hlmann, \citeyear{MeinshausenB06}; Zhao and
Yu, \citeyear{ZhaoY2006}; Zou, \citeyear{Zou2006}). Nardi and Rinaldo (\citeyear{NardiR}) considered selection
consistency of the group LASSO under an irrepresentable condition and
the bounds on the prediction and estimation errors under a restricted
eigenvalue condition (Bickel, Ritov and Tsybokov, \citeyear{BRT}; Koltchinskii,
\citeyear{KV}), assuming that the Gram matrices $X_j'X_j/n$ are proportional to
the identity matrix. Wei and Huang (\citeyear{WeiH08}) considered the sparsity and
$\ell_2$ bounds on the estimation and prediction errors of the group
LASSO under the sparse Riesz condition (Zhang and Huang, \citeyear{ZhangH2008}). They
also studied the selection property of the adaptive group LASSO using
the group LASSO as the initial estimate. The adaptive group LASSO can
be formulated in a way similar to the standard adaptive LASSO (Zou,
\citeyear{Zou2006}). Recently, there has been considerable progress in the studies of
the LASSO based on sharper versions of the restricted eigenvalue
condition (van de Geer and B\"{u}hlmann, \citeyear{GeerB2009}; Zhang, \citeyear{Zhang2009}; Ye and
Zhang, \citeyear{YeZ2010}). It would be interesting to extend these results to the
group LASSO.\vadjust{\goodbreak}

A natural question about the group LASSO is under what conditions it
will perform better than the standard LASSO. This question was
addressed by Huang and Zhang (\citeyear{HZ}), who introduced the concept of
strong group sparsity. They showed that the group LASSO is superior to
the standard LASSO under the strong group sparsity and certain other
conditions, including a group sparse eigenvalue condition. More
recently, Lounici et al. (\citeyear{LPTG2011}) conducted a detailed analysis of the
group LASSO. They established oracle inequalities for the prediction
and $\ell_2$ estimation errors of group LASSO under a restricted
eigenvalue condition on the design matrix. They also showed that the
rate of convergence of their upper bounds is optimal in a minimax
sense, up to a logarithmic factor, for all estimators over a class of
group sparse vectors. Furthermore, by deriving lower bounds for the
prediction and $\ell_2$ estimation errors of the standard LASSO they
demonstrated that the group LASSO can have smaller prediction and
estimation errors than the LASSO.

While the group LASSO enjoys excellent properties in terms of
prediction and $\ell_2$ estimation errors, its selection consistency
hinges on the assumption that the design matrix satisfies the
irrepresentable condition. This condition is, in general, difficult to
satisfy, especially in $p \gg n$ models (Zhang, \citeyear{Zhang2010}). Fan and Li
(\citeyear{FanL2001}) pointed out that the standard LASSO over-shrinks large
coefficients due to the nature of $\ell_1$ penalty. As a result, the
LASSO tends to recruit unimportant variables into the model in order to
compensate for its overshrinkage of large coefficients, and
consequently, it may not be able to distinguish variables with small to
moderate coefficients from unimportant ones. This can lead to
relatively high false positive selection rates. Leng, Lin and Wahba
(\citeyear{LengLW2006}) showed that the LASSO does not achieve selection consistency if
the penalty parameter is selected by minimizing the prediction error.
The group LASSO is likely to behave similarly. In particular, the group
LASSO may also tend to select a model that is larger than the
underlying model with relatively high false positive group selection
rate. Further work is needed to better understand the properties of the
group LASSO in terms of false positive and false negative selection
rates.

\subsection{Concave 2-Norm Group Selection}\label{sec22}
The group LASSO can be constructed by applying the $\ell_1$ penalty to
the norms of the groups. Specifically, for $\rho(t;\lam) = \lam|t|$,
the group LASSO penalty can be written as \( \lam c_j
\|\bbeta_j\|_{R_j}=\rho(\|\bbeta_j\|_{R_j}; c_j\lam). \) Other penalty
functions could be used instead. Thus a more general class of group
selection methods can be based on the criterion
%
\begin{equation}\label{gM1}
\hspace*{26pt}\frac{1}{2n}\Biggl\|\by-\sum_{j=1}^J X_j \bbeta_j\Biggr\|_2^2 + \sum_{j=1}^J
\rho
(\|\bbeta_j\|_{R_j}; c_j\lam, \gam),
\end{equation}
where $\rho(t; c_j\lam,
\gam)$ is concave in $t$. Here $\gam$ is an additional tuning parameter
that may be used to modify $\rho$. As in the definition of the group
LASSO, we assume without loss of generality that each $X_j$ is
orthonormal with $X_j'X_j/n=I_{d_j}$ and
$\|\bbeta_j\|_{R_j}=\|\bbeta_j\|_2$.

It is reasonable to use penalty functions that work well for individual
variable selection. Some possible choices include: (a) the bridge
penalty with $\rho(x;\allowbreak\lam, \gam)= \lam|x|^{\gam}, 0 < \gam\le1$
(Frank and Friedman, \citeyear{FrankF93}); (b) the SCAD penalty with $\rho(x; \lam,
\gam)=\lam\int_0^{|x|} \min\{1,\allowbreak(\gam-t/\lam)_+/(\gam-1)\}\,dt$,
$\gam>
2$ (Fan and Li, \citeyear{FanL2001}; Fan and Peng, \citeyear{FanP2004}), where for any $a \in\real$,
$a_+$ denotes its positive part, that is, $a_+=a 1_{\{a\ge0\}}$; (c)
the minimax concave penalty (MCP) with $\rho(x; \lam,
\gam)=\lam\int_0^{|x|}(1-t/(\gam\lam))_+\,dt, \gam> 1$ (Zhang, \citeyear{Zhang2010}).
All these penalties have the oracle property for individual variables,
meaning that the corresponding penalized estimators are equal to the
least squares estimator assuming the model is known with high
probability under appropriate conditions. See Huang, Horowitz and Ma (\citeyear{HuangHM}) for
the bridge penalty, Fan and Li (\citeyear{FanL2001}) and Fan and Peng (\citeyear{FanP2004}) for the
SCAD penalty and Zhang (2010) for the MC penalty. By applying these
penalties to (\ref{gM1}), we obtain the 2-norm group bridge, 2-norm
group SCAD and 2-norm group MCP, respectively. Another interesting
concave penalty is the capped-$\ell_1$ penalty $\rho(t;\lam, \gam) =
\min(\gam\lam^2/2, \lam|t|)$ with $\gam> 1$ (Zhang, \citeyear{ZhangT2011}; Shen, Zhu and Pan, \citeyear{ShenZP2011}).
However, this penalty has not been applied to the group selection
problems.

For $c_j = \sqrt{d_j}$, the group MCP and capped-$\ell_1$ pen\-alty
satisfy the invariance property
%
\begin{equation}\label{InvA}
\hspace*{25pt}\rho\bigl(\|\bbeta_j\|_2; \sqrt{d_j}\lam,
\gam\bigr)=\rho\bigl(\sqrt{d_j}\|\bbeta_j\|_2;\lam,d_j\gam\bigr).
\end{equation}
Thus the rescaling of $\lam$ can also be interpreted based on the
expression on the right-hand side of (\ref{InvA}). The multiplier
$\sqrt{d_j}$ of $\|\bbeta_j\|_2$ standardizes the group size. This
ensures that smaller groups will not be overwhelmed by larger groups.
The multiplier $d_j$ for $\gam$ makes the amount of regularization per
group proportional to its size. Thus the interpretation of $\gam$
remains the same as that in the case where group sizes are equal
to\vadjust{\goodbreak}
one. Because the MCP is equivalent to the $\ell_1$ penalty when
$\gam=\infty$, the $\ell_1$ penalty also satisfies (\ref{InvA}).
However, many other penalties, including the SCAD and $\ell_q$
penalties with $q\neq1$, do not satisfy (\ref{InvA}).

An interesting question that has not received adequate attention is how
to determine the value of $\gam$. In linear regression models with
standardized predictors, Fan and Li (\citeyear{FanL2001}) suggested using $\gam
\approx3.7$ in the SCAD penalty, and Zhang (\citeyear{Zhang2010}) suggested using
$\gam\approx2.7$ in the MCP.
Note, however, that when $\gam\rightarrow\infty$, the group MCP
converges to the group\break LASSO, and when $\gam\rightarrow1$, it
converges to the group hard threshold penalty (Antoniadis, \citeyear{Anto1996})
%
\[
\rho(t;\lam) = \lam^2-\tfrac{1}{2}(|t|-\lam)^2 1_{\{|t|\le\lam\}}.
\]
Clearly, the choice of $\gam$ has a big impact on the estimate. See
\citet{MazumderFH} and Breheny and Huang (\citeyear{BrehenyHuang2}) for further
discussion on the choice of~$\gam$.

\begin{figure*}

\includegraphics{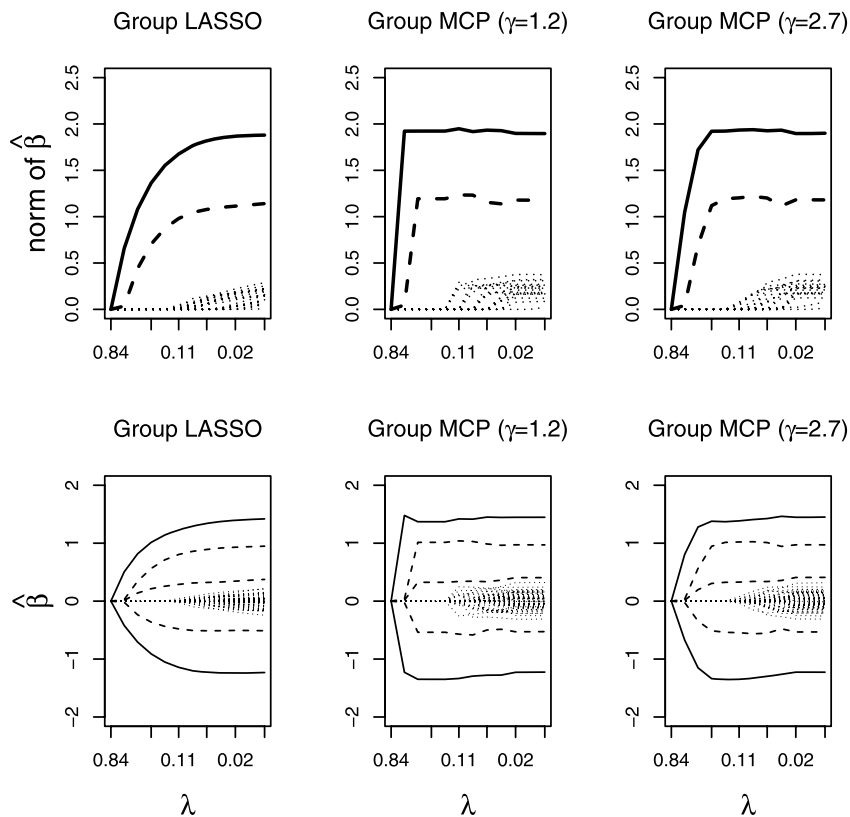}

\caption{The solution paths of the 2-norm group MCP for
$\gam=1.2, 2.7$ and $\infty$, where $\gam=\infty$ corresponds to the
group LASSO. The top panel shows the paths of the $\ell_2$ norms of
$\bbeta_j$; the bottom shows the paths of the individual coefficients.
The solid lines and dashed lines in the plots indicate the paths of the
coefficients in the nonzero groups 1 and 2, respectively. The dotted
lines represent the zero groups.}\label{fig1}\vspace*{6pt}
\end{figure*}

To illustrate this point in the grouped variable setting,
we consider a simple example 
with $J=20$ groups, in which only the first two groups have non\-zero
coefficients with $\bbeta_{1}=(-\sqrt{2}, \sqrt{2})',
\bbeta_{2}=(0.5,1,\break -0.5)'$, so $\|\bbeta_1\|_2=2$ and
$\|\bbeta_2\|_2\approx1.22$. The sizes of the groups with zero
coefficients are 3. The top panel in Figure \ref{fig1} shows the paths
of the estimated norms $\|\hbbeta_1\|$ and $\|\hbbeta_2\|$ for
$\gamma=1.2, 2.5$ and $\infty$, where $\gam=\infty$ corresponds to the
group LASSO. The bottom panel shows the solution paths of the
individual coefficients. It can be seen that the characteristics of the
solution paths are quite different for different values of $\gam$. For
the 2-norm group MCP with $\gam=1.2$ or $2.5$, there is a region in the
paths where the estimates are close to the true parameter values.
However, for the group LASSO ($\gam=\infty$), the estimates are always
biased toward zero except when $\lam=0$.

\subsection{Orthogonal Groups}
To have some understanding of the basic characteristics of the group
LASSO and nonconvex group selection methods, we consider the special
case where the groups are orthonormal with $X_j'X_k = 0, j\neq k$ and
$X_j'X_j/n = I_{d_j}$. In this case, the problem simplifies to that of
estimation in $J$ single-group models of the form \( \by=X_j\btheta+
\bveps. \) Let $\bz=X_j'\by/n$ be the least squares estimator of
$\btheta$. Without loss of generality, let $c_j=1$ below in this
section. We have $n^{-1}\|\by-X_j\btheta\|_2^2 =
\|\bz-\btheta\|_2^2+n^{-1} \|\by\|_2^2-\|\bz\|_2^2$ since
$X_j'X_j/n=I_{d_j}$. Thus the penalized least squares criterion is \(
2^{-1}\|\bz-\btheta\|_2^2+ \rho(\|\btheta\|_2;\lam, \gam). \) Denote
%
\begin{equation}\label{gLu1}
S(\bz;t) = \biggl(1-\frac{t}{\|\bz\|_2} \biggr)_{+} \bz .
\end{equation}
This expression is used in Yuan and Lin (\citeyear{YuanL2006}) for computing the group
LASSO solutions via a group coordinate descent algorithm. It is a
multivariate version of the soft-threshold operator (Donoho and
Johnstone, \citeyear{DonohoJ1994}) in which the soft-thresholding is applied to the
length of the vector, while leaving its direction unchanged. By taking
$\rho$ to be the $\ell_1$, MCP and SCAD penalties, it can be verified
that the group LASSO, group MCP and group SCAD solutions in a single
group model have the following expressions:
\begin{itemize}
\item Group LASSO:
%
\begin{equation}\label{gLu}
\hbtheta_{\mathrm{gLASSO}}(\bz; \lam)=S(\bz,
\lam).
\end{equation}
\item2-norm group MCP: for $\gamma> 1$,
\begin{eqnarray}\label{gMu}
&&\hbtheta_{\mathrm{gMCP}}(\bz;\lam, \gam)\nonumber\\ [-8pt]\\ [-8pt]
&&\quad=
\cases{
\frac{\gam}{\gam-1} S(\bz,\lam), & $\mbox{if } \|\bz\|_2 \le
\gam\lam$, \cr
\bz, & $\mbox{if } \|\bz\|_2 > \gam\lam$.
}\nonumber
\end{eqnarray}
\item2-norm group SCAD: for $\gam> 2$,
\begin{eqnarray}\label{gSu}
&&\hspace*{8pt}\hbtheta_{\mathrm{gSCAD}}(\bz;\lam, \gam)\nonumber\\ [-5pt]\\ [-8pt]
&&\hspace*{8pt}\quad=
\cases{
 S(\bz; \lam), &
$\mbox{if } \|\bz\|_2 \le2 \lam$,\vspace{+2pt}\cr
\frac{\gam-1}{\gam-2} S\bigl(\bz; \frac{\gam\lam}{\gam-1}\bigr), &
$\mbox{if } 2 \lam< \|\bz\|_2 \le\gam\lam$,\cr
\bz,&
$\mbox{if } \|\bz\|_2 > \gam\lam$.}\nonumber
\end{eqnarray}
\end{itemize}

The group LASSO solution here is simply the multivariate soft-threshold
operator. For the 2-norm\break group MCP solution, in the region $\|\bz\|_2 >
\gam\lam$, it is equal to the unbiased estimator $\bz$, and in the
remaining region, it is a scaled-up soft threshold operator. The 2-norm
group SCAD is similar to the 2-norm group MCP in that it is equal to
the unbiased estimator $\bz$ in the region $\|\bz\|_2 > \gam\lam$. In
the region $\|\bz\|_2 \le\gam\lam$, the 2-norm group SCAD is also
related to the soft threshold operator, but takes a more complicated
form than the 2-norm group MCP.

For the 2-norm\vspace*{-2pt} group MCP, $\htheta_{\mathrm{gMCP}}(\cdot;\lam,\gam)
\rightarrow\break\htheta_{\mathrm{gLASSO}}(\cdot;\lam)$ as $\gam
\rightarrow
\infty$ and $\htheta_{\mathrm{gMCP}}(\cdot; \lam, \gam) \rightarrow
H(\cdot;\lam)$ as\vadjust{\goodbreak} $\gam\to1$ for any given $\lam> 0$, where
$H(\cdot;\lam)$ is the hard-threshold operator defined as
%
\begin{equation}\label{gHu}
H(\bz; \lam)\equiv
\cases{
0, & $\mbox{if } \|\bz\|_2\le\lam$,\cr
\bz, & $\mbox{if } \|\bz\|_2> \lam$.
}
\end{equation}
Therefore, for a given $\lam> 0$, $\{\htheta_{\mathrm{gMCP}}(\cdot
;\lam,
\gam)\dvtx 1< \gam\le\infty\}$ is a family of threshold operators with
the multivariate hard and soft threshold operators at the extremes
$\gam=1$ and $\infty$.

For the 2-norm group SCAD, we have
$\htheta_{\mathrm{gSCAD}}(\cdot;\lam,\allowbreak\gam) \rightarrow
\htheta_{\mathrm{gLASSO}}(\cdot;\lam)$ as $\gam\rightarrow\infty$ and
$\htheta_{\mathrm{gSCAD}}(\cdot; \lam, \gam) \rightarrow H^*(\cdot
;\lam)$
as $\gam\to2$, where
%
\begin{equation}
\label{gHus} H^*(\bz; \lam)\equiv
\cases{
S(\bz; \lam), &$\mbox{if } \|\bz\|_2 \le2 \lam$, \cr
\bz, &$\mbox{if } \|\bz\|_2 > 2 \lam$.
}
\end{equation}
This is different from the hard threshold operator (\ref{gHu}). For a
given $\lam> 0$, $\{\htheta_{\mathrm{gSCAD}}(\cdot;\lam, \gam)\dvtx 2<
\gam
\le\infty\}$ is a family of threshold operators with $H^*$ and soft
threshold operators at the extremes $\gam=2$ and $\infty$. Note that
the hard threshold operator is not included in the group SCAD family.

The closed-form expressions given above illustrate some important
differences of the three group selection methods. They also provide
building blocks of the group coordinate descent algorithm for computing
these solutions described below.

\subsection{Computation via Group Coordinate Descent}\label{sec24}
Group coordinate descent (GCD) is an efficient approach for fitting
models with grouped penalties. The first algorithm of this kind was
proposed by Yuan and Lin (\citeyear{YuanL2006}) as a way to compute the solutions to
the group LASSO. Because the solution paths of the group LASSO are not
piecewise linear, they cannot be computed using the LARS algorithm
(Efron et al., \citeyear{EfronHJT04}).

Coordinate descent algorithms (Fu, \citeyear{Fu1998}; Friedman et al., \citeyear{FHHT2007}; Wu and
Lange, \citeyear{WuL2007}) have become widely used in the field of penalized
regression. These algorithms were originally proposed for optimization
in problems with convex penalties such as the LASSO, but have also been
used in calculating SCAD and MCP estimates (Breheny and Huang, \citeyear{BrehenyHuang2}).
We discuss here the idea behind the algorithm and its extension to the
grouped variable case.

Coordinate descent algorithms optimize an objective function with
respect to a single parameter at a time, iteratively cycling through
the parameters until convergence is reached; similarly, group
coordinate descent algorithms optimize the target function with respect
to a single group at a time, and cycles through\vadjust{\goodbreak} the groups until
convergence. These algorithms are particularly suitable for fitting
group LASSO, group SCAD and group MCP models, since all three have
simple closed-form expressions for a single-group model
(\ref{gLu})--(\ref{gSu}).

A group coordinate descent step consists of partially optimizing the
penalized least squares criterion (\ref{gL1}) or (\ref{gM1}) with
respect to the coefficients in group $j$. Define
\begin{eqnarray*}
L_j(\bbeta_j; \lam,\gam) &=& \frac{1}{2n}\Biggl\|\by-\sum_{k\neq
j}X_{k}\tbbeta_k-X_{j}\bbeta_j\Biggr\|_2^2\\
&&{} + \rho(\|\bbeta_j\|_2; c_j\lam,
\gam),
\end{eqnarray*}
where $\tbbeta$ denotes the most recently updated value of~$\bbeta$.
Denote $\tby_j=\sum_{k\neq j}X_k\tbbeta_k$ and
$\tbz_j=X_j'(\by-\tby_j)/n$. Note that $\tby_j$ represents the fitted
values excluding the contribution from group $j$, and $\tbz_j$
represents the corresponding partial residuals. Just as in ordinary
least squares regression, the value $\bbeta_j$ that optimizes
$L_j(\bbeta_j;\lam,\gam)$ is equal to the value we obtain from
regressing $\bbeta_j$ on the partial residuals. In other words, the
minimizer of $L_j(\bbeta_j;\lam,\gam)$ is given by $F(\tbz_j;
\lam,\gam)$, where $F$ is one of the solutions in (\ref{gLu}) to
(\ref{gSu}), depending on the penalty used.

Let $\tbbeta^{(0)}=(\tbbeta_1^{(0)\prime},
\ldots,\tbbeta_J^{(0)\prime})'$ be the initial value, and let $s$
denote the iteration. The GCD algorithm consists of the following
steps:

\begin{itemize}
\item[] Step 1. Set $s=0$. Initialize vector of residuals
$\bfr=\by-\tby$, where $\tby=\sum_{j=1}^J X_{j}\tbbeta_j^{(0)}$.

\item[] Step 2. For $j =1, \ldots, J$, carry out the following
calculations:
\begin{itemize}
\item[(a)] calculate \( \tbz_j = n^{-1}X_j'\bfr+\tbbeta_j^{(s)}; \)
\item[(b)] update $\tbbeta_j^{(s+1)}=F(\tbz_j; \lam,\gam)$,

\item[(c)] update $ r \gets r -
X_j(\tbbeta_j^{(s+1)}-\tbbeta_j^{(s)})$.
\end{itemize}

\item[] Step 3. Update $s \gets s+1$.

\item[] Step 4. Repeat steps 2--3 until convergence.
\end{itemize}

The update in Step 2(c) ensures that $\bfr$ always holds the current
values of the residuals, and is therefore ready for Step 2(a) of the
next cycle. By taking $F(\cdot; \lam, \gam)$ to be
$\hbtheta_{\mathrm{gLASSO}}(\cdot; \lam)$, $\hbtheta_{\mathrm
{gMCP}}(\cdot;
\lam, \gam)$ and $\hbtheta_{\mathrm{gSCAD}}(\cdot; \lam, \gam)$ in
(\ref{gLu}) to (\ref{gSu}), we obtain the solutions to the group LASSO,
group MCP and group SCAD, respectively. The algorithm has two
attractive features. First, each step is very fast, as it involves
only relatively simple calculations. Second, the algorithm is stable,
as each step is guaranteed to decrease the objective function (or leave
it unchanged).

The above algorithm computes $\hbbeta$ for a given $(\lam, \gam)$ pair;
to obtain pathwise solutions, we can use the algorithm repeatedly over
a grid of $(\lam, \gam)$ values. For~a~given value of $\gam$, we can
start at $\lam_{\max}=\break\max_j\{\|n^{-1}X_j\by\|_2/c_j\}$, for which
$\hbbeta$ has the solution $0$, and proceed along the grid using the
value of $\hbbeta$ at the previous point in the $\lam$-grid as the
initial value for the current point in the algorithm. An alternative
approach is to use the group LASSO solution (corresponding to
$\gam=\infty$) as the initial value as we decrease $\gam$ for each
value of $\lam$. See \citet{MazumderFH} for a detailed description
of the latter approach in the nongrouped case.

The results of Tseng (\citeyear{Tseng2001}) establish that the algorithm converges to a
minimum. For the group LASSO, which has a convex objective function,
the algorithm therefore converges to the global minimum. For group
SCAD and group MCP, convergence to a local minimum is possible. See
also Theorem 4 of \citet{MazumderFH} for the nongrouped case.

The availability of the explicit expression in step 2(b) of the
algorithm depends on the choice of $R_j=X_j'X_j/n$ in (\ref{gL1}) or
(\ref{gM1}). If a different norm is used, then the groups are not
orthonormal, and there are no explicit solutions to the problem.
Without closed-form solutions, step 2(b) must be solved using numerical
optimization. Algorithms proposed for computing the group LASSO
solutions without using $R_j=X_j'X_j/n$ include Friedman et al. (\citeyear{FHHT2007}),
\citet{JOV} and Liu and Ye (\citeyear{LY2010}). For generalized linear
models, the group coordinate descent can be applied based on quadratic
approximations to the log-likelihood in the objective function (\cite{MGB2008}).

\section{Bi-Level Selection}
The methods described in Section \ref{sec2} produce estimates that are sparse at
the group level and not at the level of individual variables. Within a
group, there are only two possibilities for the selection results based
on these methods: either all of the variables are selected, or none of
them are. This is not always appropriate for the data.

For example, consider a genetic association study in which the
predictors are indicators for the presence of genetic variation at
different markers. If a genetic variant located in a gene is
associated with the disease, then it is more likely that other variants
located in the same gene will also be associated with the disease---the\vadjust{\goodbreak}
predictors have a grouping structure. However, it is not
necessarily the case that \emph{all} variants within that gene are
associated with the disease. In such a study, the goal is to identify
important \emph{individual} variants, but to increase the power of the
search by incorporating grouping information.\looseness=1

In this section, we discuss bi-level selection methods, which are
capable of selecting important groups as well as important individual
variables within those groups. The underlying assumption is that the
model is sparse at both the group and individual variable levels. That
is, the nonzero group coefficients $\bbeta_j$ are also sparse. It
should be noted, however, that less work has been done on bi-level
selection than on group LASSO, and there are still many unanswered
questions.

\subsection{Concave 1-Norm Group Penalties}
As one might suspect, based on analogy with\break LASSO and ridge regression,
it is possible to construct penalties for bi-level selection by
starting with the $\ell_1$ norm instead of the $\ell_2$ norm. This
substitution is not trivial, however: a na\"ive application of the
LASSO penalty to the $\ell_1$ norm of a group results in the original
LASSO, which obviously has no grouping properties.

Applying a concave penalty to the $\ell_1$ norm of a group, however,
does produce an estimator with grouping properties, as suggested by
Huang et al. (\citeyear{HuangMXZ}), who proposed the group bridge penalty. The 1-norm
group bridge applies a bridge penalty to the $\ell_1$ norm of a group,
resulting in the criterion
%
\begin{equation}\label{GBA1}
\frac{1}{2n}\Biggl\|\by-\sum_{j=1}^J
X_j\bbeta_j\Biggr\|_2^2 + \lam\sum_{j=1}^J c_j\|\bbeta_{_j}\|_1^{\gam},
\end{equation}
where $\lam>0$ is the regularization parameter, $\gam\in(0,1)$ is the
bridge index and $\{c_j\}$ are constants that adjust for the dimension
of group $j$. For models with standardized variables, a reasonable
choice is $c_j= |d_j|^{\gamma}$. When $d_j=1, 1\le j \le J$,
(\ref{GBA1}) simplifies to the standard bridge criterion. The method
proposed by Zhou and Zhu (\citeyear{ZhouZ2010}) can be considered a special case of
group bridge with $\gam=0.5$. A general composite absolute penalty
based on $\ell_{q}$ norms was proposed by Zhao, Rocha and Yu (\citeyear{ZhaoRY2009}).

Huang et al. (\citeyear{HuangMXZ}) showed that the global group bridge solution is
group selection consistent under certain regularity conditions. Their
results allow $p \to\infty$ as $n \to\infty$ but require $p < n$. In
contrast to the group LASSO, the selection consistency of group bridge
does not require an irrepresentable-type condition. However, no results
are available for the group bridge in the $J \gg n$ settings.

In principle, we could apply other concave penalties to the group
$\ell_1$ norm as well, leading to the more general penalized criterion
%
\begin{equation}\label{gM2}
\hspace*{25pt}\frac{1}{2n}\Biggl\|\by-\sum_{j=1}^J X_j
\bbeta_j\Biggr\|_2^2 +
\sum_{j=1}^J \rho(\|\bbeta_j\|_1; c_j\lam, \gam).
\end{equation}
Choosing
$\rho$ to be the SCAD or MCP penalty in (\ref{gM2}) would seem
particularly promising, but to our knowledge, these estimators have not
been studied.

\begin{figure*}

\includegraphics{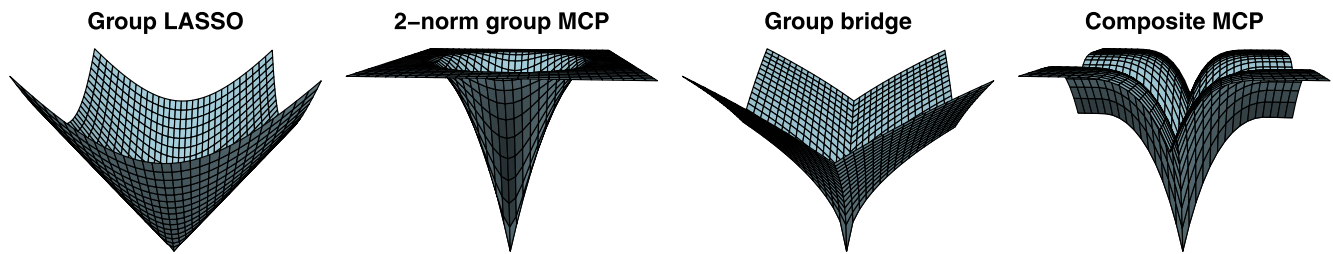}

\caption{The group LASSO, group bridge and composite mcp penalties for a
two-predictor group. Note that where the penalty comes to a point or
edge, there is the possibility that the solution will take on a sparse
value; all penalties come to a point at $\mathbf{0}$, encouraging
group-level sparsity, but only group bridge and composite MCP allow for
bi-level selection.}\label{fig:penshapes}
\end{figure*}

\subsection{Composite Penalties}
An alternative way of thinking about concave\break \mbox{1-norm} group penalties is
that they represent the composition of two penalties: a concave
group-level penalty and an individual variable-level 1-norm\break penalty.
It is natural, then, to also consider the composition of concave
group-level penalties with other individual variable-level penalties.
This framework was proposed in Breheny and Huang (\citeyear{BrehenyHuang1}), who described
grouped penalties as consisting of an outer penalty $\rho_O$ applied to
a sum of inner penalties $\rho_I$. The penalty applied to a group of
predictors is therefore written as
%
\begin{equation}\label{fmw1}
\rho_O \Biggl(\sum_{k=1}^{d_j} \rho_I(\vert\beta
_{jk}\vert) \Biggr),
\end{equation}
where $\beta_{jk}$ is the $k$th member of the $j$th group,
and the partial derivative with respect to the $jk$th covariate is
%
\begin{equation}\label{fmw2}
\rho_O' \Biggl(\sum_{k=1}^{d_j} \rho_I(\vert\beta_{jk}\vert) \Biggr)
\rho_I'(\vert\beta_{jk}\vert) .
\end{equation}
Note that the group bridge fits into this framework with an outer
bridge penalty and an inner LASSO penalty, as does the group LASSO with
an outer bridge penalty and an inner ridge penalty.

From (\ref{fmw1}), we can view group penalization as applying a rate of
penalization to a predictor that consists of two terms: the first
carries information regarding the group; the second carries information
about the individual predictor. Whether or not a variable enters the
model is affected both by its individual signal and by the collective
signal of the group that it belongs to. Thus, a variable with a
moderate individual signal may be included in a model if it belongs to
a group containing other members with strong signals, or may be
excluded if the rest of its group displays little association with the
outcome.

An interesting special case of the composite pen\-alty is using the MCP
as both the outer and inner penalties, which we refer to as the
composite MCP (this penalty was referred to as ``group MCP'' in Breheny
and Huang (\citeyear{BrehenyHuang1}); we use ``composite MCP'' both to better reflect the
framework and avoid confusion with the 2-norm group MCP of Section~\ref{sec22}).

The composite MCP uses the criterion
\begin{eqnarray}\label{gm}
&&\frac{1}{2n}
\Biggl\|\by-\sum_{j=1}^J X_j \bbeta_j\Biggr\|_2^2\nonumber\\ [-8pt]\\ [-8pt]
&&\quad{} + \sum_{j=1}^J
\rho_{\lambda,\gamma_O} \Biggl( \sum_{k=1}^{d_j} \rho_{\lambda,\gamma_I}
(\vert\beta_{jk}\vert) \Biggr),\nonumber
\end{eqnarray}
where $\rho$ is the MCP
penalty and
$\gamma_O$, the tuning parameter of the outer penalty, is chosen to be
$d_j\gamma_I\lambda/2$ in order to ensure that the group level penalty
attains its maximum if and only if each of its components are at their
maximum. In other words, the derivative of the outer penalty reaches 0
if and only if $\vert\beta_{jk}\vert\geq\gamma_I\lambda\
\forall k \in\{1,\ldots, d_j\}$.

Figure  \ref{fig:penshapes} shows the group LASSO, 2-norm group MCP,
1-norm group Bridge and composite MCP penalties for a two-predictor
group. Note that where the penalty comes to a point or edge, there is
the possibility that the solution will take on a sparse value; all
penalties come to a point at $\mathbf{0}$, encouraging group-level
sparsity, but only group bridge and composite MCP allow for bi-level
selection. In addition, one can see that the MCP penalties are capped,
while the group LASSO and group bridge penalties are not. Furthermore,
note that the individual variable-level penalty for the composite MCP
is capped at a level below that of the group; this limits the extent to
which one variable can dominate the penalty of the entire group. The
2-norm group MCP does not have this property. This illustrates the two
rationales of composite MCP: (1) to avoid overshrinkage by allowing
covariates to grow large, and (2) to allow groups to remain sparse
internally. The 1-norm group bridge allows the presence of a single
large predictor to continually lower the entry threshold of the other
variables in its group. This property, whereby a single strong
predictor draws others into the model, prevents the group bridge from
achieving consistency for the selection of individual variables.

\begin{figure*}[b]

\includegraphics{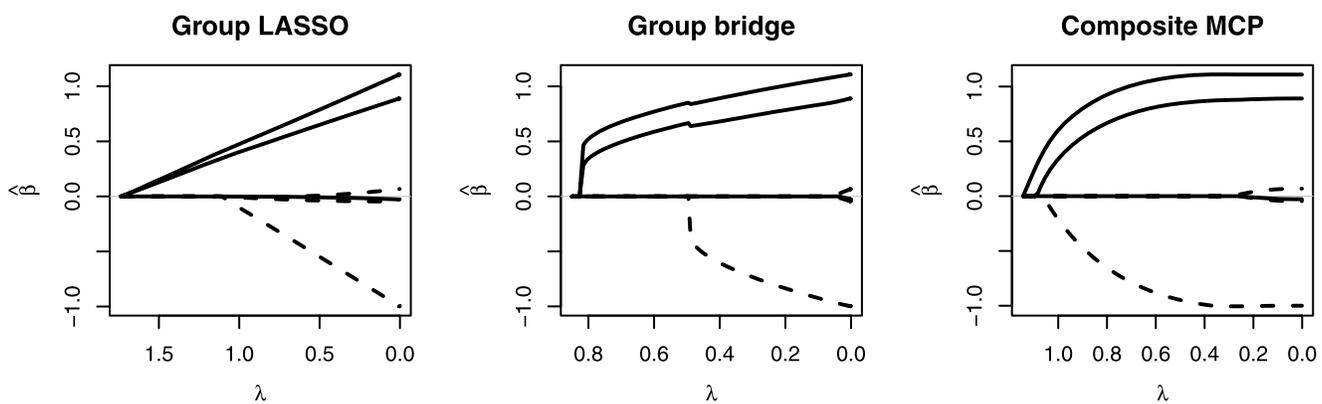}

\caption{Coefficient paths from 0 to $\lambda
_{\max}$ for group LASSO, 2-norm group MCP, 1-norm group bridge,
and composite MCP for a simulated data set featuring two groups, each
with three covariates. In the underlying data-generating mechanism, the
group represented by solid lines has two covariates with coefficients
equal to 1 and the other equal to 0; the group
represented by dashed lines has two coefficients equal to 0 and the
other equal to $-1$.}\label{fig:paths}
\end{figure*}

Figure \ref{fig:paths} shows the coefficient paths from
$\lambda_{\max}$ down to 0 for group LASSO, 1-norm group bridge,
and composite MCP for a simulated data set featuring two groups, each
with three covariates. In the underlying model, the group represented
by solid lines has two covariates with coefficients equal to 1 and the
other equal to 0; the group represented by dashed lines has two
coefficients equal to 0 and the other equal to $-1$. The figure reveals
much about the behavior of grouped penalties. In particular, we note
the following: (1) Even though each of the nonzero coefficients is of
the same magnitude, the coefficients from the more significant solid
group enter the model more easily than the lone nonzero coefficient
from the dashed group. (2) This phenomenon is less pronounced for
composite MCP, which makes weaker assumptions about grouping. (3) For
composite MCP at $\lambda\approx0.3$, all of the variables with true
zero coefficients have been eliminated while the remaining coefficients
are unpenalized. In this region, the composite MCP approach is
performing as well as the oracle model. (4) In general, the
coefficient paths for these group penalization methods are continuous,
but are not piecewise linear, unlike those for the LASSO.

Although composite penalties do not, in general, have closed-form
solutions in single-group models like the penalties in Section \ref{sec2}, the
idea of group coordinate descent can still be used. The main
complication is in step 2(b) for the algorithm described in Section
\ref{sec24}, where the single-group solutions need to be solved numerically.
Another approach is based on a local coordinate descent algorithm
(Breheny and Huang, \citeyear{BrehenyHuang1}). This algorithm first uses a local linear
approximation to the penalty function (Zou and Li, \citeyear{ZouL2008}). After
applying this approximation, in any given coordinate direction the
optimization problem is equivalent to the one-dimensional LASSO, which
has the soft-threshold operator as its solution. The thresholding
parameter $\lambda$ in each update is given by expression (\ref{fmw2}).
Because the penalties involved are concave on $[0,\infty)$, the linear
approximation is a majorizing function, and the algorithm thus enjoys
the descent property of MM algorithms (\cite{Lange2000}) whereby the
objective function is guaranteed to decrease at every iteration.
Further details may be found in Breheny and Huang (\citeyear{BrehenyHuang1}). These
algorithms have been implemented in the \texttt{R} package
\texttt{grpreg}, available at \url{http://cran.r-project.org}. The
package com-\break putes~the group LASSO, group bridge and composite MCP
solutions for linear regression and logistic regression models.

\subsection{Additive Penalties}
Another approach to achieving bi-level selection is to add an $\ell_1$
penalty to the group LASSO (Wu and Lange, \citeyear{WuL2007}; Friedman, Hastie and
Tibshirani, \citeyear{FHT2010}).
%
\begin{equation}\label{sgL}
\hspace*{26pt}\frac{1}{2n} \Biggl\|\by-\sum_{j=1}^J
X_j\bbeta_j\Biggr\|_2^2 + \lam_1\|\bbeta\|_1 + \lam_2 \sum_{j=1}^J
\|\bbeta_j\|_2,
\end{equation}
where $\lam_1\ge0$ and $\lam_2 \ge0$ are
regularization parameters. The above objective function has the
benefit of being convex, eliminating the possibility of convergence to
a local minimum during model fitting. The group coordinate descent
algorithm can no longer be applied, however, as the orthonormalization
procedure described in Section \ref{sec2} will not preserve the sparsity
achieved by the $\ell_1$ penalty once the solution is transformed back
to the original variables. \citet{PWH2009}, \citet{FHT2010}
and Zhou et al. (\citeyear{ZhouSSL2010}) have proposed algorithms for solving this
problem without requiring orthonormalization.

In principle, the group LASSO portion of the pen\-alty could be replaced
with any of the convex 2-norm group penalties of Section \ref{sec22}; likewise
the $\ell_1$ penalty could be replaced by, say, MCP or SCAD. These
possibilities, to the best of our knowledge, have not been explored.
Further work is needed to study the properties of this class of
estimators and compare their performance with other methods.

\subsection{Example: Genetic Association}
We now give an example from a genetic association study where bi-level
selection is an important goal of the study.
The example involves data from a case-control study of age-related
macular degeneration consisting of 400 cases and 400 controls, and was
analyzed in Breheny and Huang (\citeyear{BrehenyHuang1}). The analysis is confined\vadjust{\goodbreak} to 30
genes containing 532 markers that previous biological studies have
suggested may be related to the disease.

We analyze the data with the group LASSO,\break 1-norm group bridge and
composite MCP methods by considering markers to be grouped by the gene
they belong to. Penalized logistic regression models were fit assuming
an additive effect for all markers (homozygous dominant~$=$ 2,
heterozygous~$=$ 1, homozygous recessive~$=$ 0). In addition to the group
penalization methods, we analyzed these data using a traditional
one-at-a-time approach (single-marker analysis), in which univariate
logistic regression models were fit and marker effects screened using a
$p < 0.05$ cutoff. Ten-fold cross-validation was used to select
$\lambda$, and to assess accuracy (for the one-at-a-time approach,
predictions were made from an unpenalized logistic regression model fit
to the training data using all the markers selected by individual
testing). The results are presented in Table \ref{tab:app}.

\begin{table}
\caption{Application of the three group penalization
methods and a~one-at-a-time method to a genetic association data set.
CV~error is the average number of misclassification errors~over the ten
validation folds}\label{tab:app}
\begin{tabular*}{\columnwidth}{@{\extracolsep{\fill}}lccc@{}}
\hline
& \textbf{Genes} & \textbf{Markers} & \textbf{Cross-validation}\\
& \textbf{selected} & \textbf{selected} & \textbf{error} \\
\hline
One-at-a-time & 19 & \phantom{0}49 & 0.441 \\
Group LASSO & 17 & 435 & 0.390 \\
Group bridge & \phantom{0}3 & \phantom{0}20 & 0.400 \\
Composite MCP & \phantom{0}8 & \phantom{0}11 & 0.391 \\
\hline
\end{tabular*}
\end{table}

Table \ref{tab:app} suggests the benefits of using group penalization
regression approaches as opposed to one-at-a-time approaches: the three
group penalization methods achieve lower test error rates and do so
while selecting fewer genes (groups). Although the error rates of
$\approx$40\% indicate that these 30 genes likely do not include SNPs
that exert a large effect on an individual's chances of developing
age-related macular degeneration, the fact that they are well below the
50\% that would be expected by random chance demonstrates that these
genes do contain SNPs related to the disease. The very different
nature of the selection properties of the three group penalization
methods are also clearly seen. Although group LASSO achieves low
misclassification error, it selects 17 genes out of 30 and 435 markers
out of 532, failing to shed light on the most important genetic
markers. The bi-level selection methods achieve comparable error rates
with a much more sparse set of predictors: group bridge identifies
3\vadjust{\goodbreak}
promising genes out of 30 candidates, and composite MCP identifies 11
promising SNPs out of 532.

\section{Oracle Property of the 2-Norm Group MCP}\label{Oracle}

In this section, we look at the selection properties of the 2-norm
group MCP estimator $\hbbeta(\lam, \gam)$, defined as the global
minimizer of (\ref{gM1}) with $c_j = \sqrt{d_j}$, when $\rho$ is taken
to be the MCP penalty. We provide sufficient conditions under which the
2-norm group MCP estimator is equal to the oracle least squares
estimator defined at (\ref{ORA}) below. Our intention is to give some
preliminary theoretical justification for this concave group selection
method under reasonable conditions, not necessarily to obtain the best
possible theoretical results or to provide a systematic treatment of
the properties of the concave group selection methods discussed in this
review.

Let $X=(X_1, \ldots, X_J)$ and $\Sigma= X'X/n$. For any $A \subseteq
\{1, \ldots, J\}$, denote
\[
X_A = (X_j, j \in A), \quad\Sigma_A = X_A'X_A/n.
\]
Let the true value of the regression coefficients be
$\bbeta^o=(\bbeta_1^{o\prime}, \ldots, \bbeta_J^{o\prime})'$. Let
$S=\{j\dvtx \|\bbeta_j^o\|_2 \neq0, 1\le j \le J\}$, which is the set of
indices of the groups with nonzero coefficients in the underlying
model. Let $\beta_*^o= \min\{\|\bbeta_j^o\|_2/\sqrt{d_j}\dvtx j\in S\}$
and set $\beta_*^o=\infty$ if $S$ is empty.
Define
%
\begin{equation}\label{ORA}
\hspace*{20pt}\hbbeta^o=\argmin\limits_b\{ \|\by-X\bb\|_2^2\dvtx \bb_j =\bzero
\ \forall\! j \notin S\}.
\end{equation}
This is the oracle least squares estimator. Of course, it is not a real
estimator, since the oracle set is unknown.

Let $d_{\max}=\max\{d_j\dvtx 1\le j \le J\}$ and $d_{\min}=\break \min\{d_j\dvtx
1\le
j \le J\}$. For any $A \subseteq\{1, \ldots, J\}$, denote
$d_{\min}(A)=\min\{d_j\dvtx j \in A\}$ and $d_{\max}(A)=\max\{d_j\dvtx\break j \in A\}$.
Here $d_{\min}(A)=\infty$ if $A$ is empty. Let $c_{\min}$ be the
smallest eigenvalue of $\Sigma$, and let $c_1$ and $c_2$ be the
smallest and largest eigenvalues of $\Sigma_{S}$, respectively.

We first consider the case where the 2-norm group MCP objective
function is convex. This necessarily requires $c_{\min} > 0$. Define
the function
%
\begin{eqnarray}\label{Defh}
h(t, k)=\exp\bigl(-k\bigl(\sqrt{2t-1}-1\bigr)^2/4 \bigr),\nonumber\\ [-8pt]\\ [-8pt]
  \eqntext{t > 1, k=1, 2,\ldots.}
\end{eqnarray}
This function arises from an upper bound for
the tail probabilities of chi-square distributions given in Lemma \ref{LemChi2} in
the \hyperref[appendix]{Appendix}, which is based on an exponential inequality for
chi-square random variables of Laurent and Massart (\citeyear{LaurentM2010}). Let
%
\begin{equation}\label{DefEta1}
\eta_{1n}(\lam) =(J-|S|) h \bigl( \lam^2 n/\sigma^2,
d_{\min}(S^c) \bigr)
\end{equation}
and
%
\begin{equation}\label{DefEta2}
\hspace*{20pt}\eta_{2n}(\lam) = |S|h \bigl(c_1
n(\beta_*^o-\gam\lam)^2/\sigma^2, d_{\min}(S) \bigr).
\end{equation}

\begin{theorem}\label{ThmA}
Suppose $\veps_1, \ldots, \veps_n$ are independent and
identically distributed as $N(0, \sigma^2)$. Then for any $(\lam,
\gam)$ satisfying $\gamma> 1/c_{\min}$, $\beta_*^o > \gamma\lam$ and
$n \lam^2 > \sigma^2$, we have
\[
\rP\bigl(\hbbeta(\lam, \gam) \neq\hbbeta^o \bigr) \le\eta_{1n}(\lam)+
\eta_{2n}(\lam).
\]
\end{theorem}

The proof of this theorem is given in the \hyperref[appendix]{Appendix}. It provides an
upper bound on the probability that $\hbbeta(\lam, \gam)$ is not equal
to the oracle least squares estimator. The condition $\gam>
1/c_{\min}$ ensures that the 2-norm group MCP criterion is strictly
convex. This implies
$\hbbeta(\lam, \gam)$ is uniquely characterized by the
Ka\-rush--Kuhn--Tucker conditions. The condition $n\lam^2 > \sigma^2$
requires that $\lam$ cannot be too small.

Let
%
\begin{eqnarray}\label{DefLamTau}
\hspace*{16pt}\lam_n &=& \sigma\bigl( 2\log(\max\{J-|S|,1\})\nonumber\\
\hspace*{16pt}&&\hspace*{53.4pt}{}/(n
d_{\min}(S^c))\bigr)^{1/2}
\quad\mbox{and}\\
\hspace*{16pt}\tau_n &=& \sigma\sqrt{2\log(\max\{|S|, 1\}) /(n c_1 d_{\min}(S))}.\nonumber
\end{eqnarray}

The following corollary is an immediate consequence of Theorem
\ref{ThmA}.

\begin{corollary}\label{ColA}
Suppose that the conditions of Theorem \ref{ThmA} are
satisfied. Also suppose that $\beta_*^o \ge\gamma\lam+ a_n\tau_n$ for
$a_n \rightarrow\infty$ as $n \rightarrow\infty$. If $\lam\ge a_n
\lam_n$, then
\[
\rP\bigl(\hbbeta(\lam,\gam) \neq\hbbeta^o\bigr) \rightarrow0 \quad\mbox{as } n
\rightarrow\infty.
\]
\end{corollary}

By Corollary \ref{ColA}, the 2-norm group MCP estimator behaves like
the oracle least squares estimator with high probability. This of
course implies it is group selection consistent. For the standard LASSO
estimator, a sufficient condition for its sign consistency is the
strong irrepresentable condition (Zhao and Yu, \citeyear{ZhaoY2006}). Here a similar
condition holds automatically due to the form of the MCP. Specifically,
let $\bbeta_S^o=(\bbeta_j^{o'}\dvtx j\in S)'$. Then an extension of the
irrepresentable condition to the present setting is, for some $0 <
\delta< 1$,
\begin{eqnarray}\label{IrrepA}
&&\max_{j\notin S}\|X_j'X_S(X_S'X_S)^{-1}
\drho(\bbeta_S^o; \lam, \gam)/\lam\|_2\nonumber\\ [-8pt]\\ [-8pt]
&&\quad \le1- \delta,\nonumber
\end{eqnarray}
where
$\drho(\bbeta_S^o;\lam, \gam)= (\drho(\|\bbeta_j^o\|_2;
\sqrt{d_j}\lam,\gam)\bbeta_j^{o\prime}/\|\bbeta_j^o\|_2\dvtx\break j \in S)'$
with
\[
\drho\bigl(\|\bbeta_j^o\|_2;\sqrt{d_j}\lam,\gam\bigr) = \lam
\bigl(1-\|\bbeta_j^o\|_2/\bigl(\sqrt{d_j}\gam\lam\bigr) \bigr)_+ .
\]
Since it is assumed that $\min_{j\in S}\|\bbeta_j^o\|_2/\sqrt{d_j} >
\gam\lam$, we have $\drho(\|\bbeta_j^o\|_2; \sqrt{d_j}\lam,\gam)=0$
for all $j \in S$. Therefore, (\ref{IrrepA}) always holds.

We now consider the high-dimensional case where $J > n$.
We require the sparse Riesz condition, or SRC (Zhang and Huang, \citeyear{ZhangH2008}),
which is a form of sparse eigenvalue condition. We say that $X$
satisfies the SRC with rank $d^*$ and spectrum bounds $\{c_*, c^*\}$ if
%
\begin{eqnarray}\label{SRCb}
0 < c_*\le\|X_A \bu\|_2^2/n \le c^*< \infty\nonumber\\ [-8pt]\\ [-8pt]
\eqntext{\forall\! A
\mbox{ with } |A| \le d^*, \|\bu\|_2=1.}
\end{eqnarray}
We refer to this condition as $\operatorname{SRC}(d^*, c_*, c^*)$.

Let $K_* = (c^*/c_*)-(1/2)$, $m_*=K_*|S|$ and $\xi=1/(4 c^*d_s)$, where
$d_s=\max\{d_{\max}(S), 1\}$. Define
\begin{eqnarray}\label{DefPi1s3}
\eta_{3n}(\lam)&=& (J-|S|)^{m_*}
\frac{e^{m_*}}{m_{*}^{m_*}}\nonumber\\ [-8pt]\\ [-8pt]
&&{}\cdot h(\xi n\lam^2\sigma^{-2}/d_{\max}, m_*
d_{\max}).\nonumber
\end{eqnarray}
Let $\eta_{1n}$ and $\eta_{2n}$ be as in (\ref{DefEta1}) and
(\ref{DefEta2}).

\begin{theorem}\label{ThmB}
Suppose $\veps_1, \ldots, \veps_n$ are independent and
identically distributed as $N(0, \sigma^2)$, and $X$ satisfies the
$\operatorname{SRC}(d^*, c_*, c^*)$ in (\ref{SRCb}) with $d^* \ge(K_*+1)|S|d_s$.
Then for any $(\lam, \gam)$ satisfying $\beta_*^o > \gam\lam$,\break
$n\lam^2 \xi> \sigma^2 d_{\max}$ and \(\gam\ge
c_*^{-1}\sqrt{4+(c_*/c^*)}\), we have
\[
\rP\bigl( \hbbeta(\lam, \gam) \neq\hbbeta^o \bigr) \le\eta_{1n}(\lam) +
\eta_{2n}(\lam)+\eta_{3n}(\lam).
\]
\end{theorem}

Letting
\[
\lam_n^* = 2\sigma\sqrt{2c^*d_s \log(J-|S|)/n}
\]
and $\tau_n$ be as in (\ref{DefLamTau}), Theorem \ref{ThmB} has the
following corollary.

\begin{corollary}
\label{ColB} Suppose the conditions of Theorem \ref{ThmB} are
satisfied. Also suppose $\beta_*^o \ge\gamma\lam+ a_n\tau_n$ for $a_n
\rightarrow\infty$ as $n \rightarrow\infty$. Then if $\lam\ge a_n
\lam_n^*$,
\[
\rP\bigl( \hbbeta(\lam, \gam) \neq\tbbeta^o \bigr) \rightarrow0 \quad \mbox{as } n
\rightarrow\infty.
\]
\end{corollary}
Theorem \ref{ThmB} and Corollary \ref{ColB} provide sufficient
conditions for the selection consistency of the global 2-norm group MCP
estimator in the $J \gg n$ situations. For example, we can have
$J-|S|=\exp\{o(n/\break (c^* d_s))\}$. The condition $ n\lam^2 \xi> \sigma^2
d_{\max}$ is stronger than the corresponding condition $n\lam^2 >
\sigma^2$ in Theorem \ref{ThmA}. The condition $\gam\ge
c_*^{-1}\sqrt{4+(c_*/c^*)}$ ensures that the group MCP criterion is
convex in any $d^*$-dimensional subspace. It is stronger than the
minimal sufficient condition $\gamma> 1/c_*$ for convexity in
$d^*$-dimensional subspaces. These reflect the difficulty and extra
efforts needed in reducing a $p$-dimensional problem to a
$d^*$-dimensional problem. The SRC in (\ref{SRCb}) guarantees that the
model is identifiable in a lower $d^*$-dimensional space.

The results presented above are concerned with the global solutions.
The properties of the local solutions, such as those produced by the
group coordinate descent algorithm, to concave 2-norm or 1-norm
penalties remain largely unknown in models with $J \gg n$. An
interesting question is under what conditions the local solutions are
equal to or sufficiently close to the global solutions so that they are
still selection consistent. In addition, the estimation and prediction
properties of these solutions have not been studied. We expect that the
methods of Zhang and Zhang (\citeyear{ZhangZ2011}) in studying the properties of concave
regularization will be helpful in group and bi-level selection
problems.

\section{Applications}
We now give a review of some applications of the group selection
methods in several statistical modeling and analysis problems,
including nonparametric additive models, semiparametric partially
linear models, seemingly unrelated regressions and multi-task learning
and genetic and genomic data analysis.

\subsection{Nonparametric Additive Models} Let $(y_i,
\bx_i), i=1, \ldots, n$ be random vectors that are independently and
identically distributed as $(y, \bx)$, where $y$ is a response
variable, and $\bx=(x_{1}, \ldots, x_{p})'$ is a $p$-dimensional
covariate vector. The nonparametric additive model (Hastie and
Tibshirani, \citeyear{HT}) posits that
%
\begin{equation}
\label{GamA} y_i= \mu+ \sum_{j=1}^{p} f_j(x_{ij}) + \veps_i,\quad 1\le
i\le
n,
\end{equation}
where $\mu$ is an intercept term, $x_{ij}$ is the $j$th component of
$x_{i}$, the $f_j$'s are unknown functions and $\veps_i$ is an
unobserved random variable with mean zero and finite variance
$\sigma^2$. Suppose that some of the additive components $f_j$ are
zero. The problem is to select the nonzero components and estimate
them. Lin and Zhang (\citeyear{LinZhang2006}) proposed the component selection and
smoothing operator (COSSO) method that can be used for selection and
estimation in (\ref{GamA}). The COSSO can be viewed as a group LASSO
procedure in a reproducing kernel Hilbert space. For fixed $p$, they
studied the rate of convergence of the COSSO estimator in the additive
model. They also showed that, in the special case of a tensor product
design, the COSSO correctly selects the non-zero additive components
with high probability. Zhang and Lin (\citeyear{ZhangL2006}) considered the COSSO for
nonparametric regression in exponential families. Meier, van de Geer
and B\"{u}hlmann (\citeyear{MerierGB}) proposed a variable selection method in
(\ref{GamA}) with $p \gg n$ that is closely related to the group LASSO.
They give conditions under which, with high probability, their
procedure selects a set of the nonparametric components whose distance
from zero in a certain metric exceeds a specified threshold under a
compatibility condition. Ravikumar et al. (\citeyear{RLLW2009}) proposed a penalized
approach for variable selection in (\ref{GamA}). In their theoretical
results on selection consistency, they assume that the eigenvalues of
a ``design matrix'' be bounded away from zero and infinity, where the
``design matrix'' is formed from the basis functions for the non\-zero
components. Another critical condition required in their paper is
similar to the irrepresentable condition of Zhao and Yu (\citeyear{ZhaoY2006}). Huang,
Horowitz and Wei (\citeyear{HHW}) studied the group LASSO and adaptive group
LASSO for variable selection in (\ref{GamA}) based on a spline
approximation to the nonparametric components. With this approximation,
each nonparametric component is represented by a linear combination of
spline basis functions. Consequently, the problem of component
selection becomes that of selecting the groups of coefficients in the
linear combinations. They provided conditions under which the group
LASSO selects a model whose number of components is comparable with the
underlying model, and the adaptive group LASSO selects the nonzero
components correctly with high probability and achieves the optimal
rate of convergence.

\subsection{Structure Estimation in Semiparametric Regression
Models}
Consider the semiparametric partially linear model (Engle
et al., \citeyear{En86})
\begin{eqnarray}
\label{PLM1} y_i&=& \mu+ \sum_{j\in S_1} \beta_j x_{ij}\nonumber\\ [-8pt]\\ [-8pt]
 &&{}+ \sum_{j\in
S_2} f_j(x_{ij})+ \veps_i,\quad 1 \le i \le n,\nonumber
\end{eqnarray}
where $S_1$ and $S_2$ are mutually exclusive and complementary subsets
of $\{1, \ldots, p\}$, $\{\beta_j\dvtx j \in S_1\} $ are regression\vadjust{\goodbreak}
coefficients of the covariates with indices in $S_1$ and $(f_{j}\dvtx j
\in S_2)$ are unknown functions. The most important assumption in the
existing methods for the estimation in partially linear models is that
$S_1$ and $S_2$ are known a priori. This assumption underlies the
construction of the estimators and investigation of their theoretical
properties in the existing methods (H\"{a}rdle, Liang and Gao, \citeyear{HLG};
Bickel et al., \citeyear{BKRW}). However, in applied work, it is rarely known in
advance which covariates have linear effects and which have nonlinear
effects. Recently, Zhang, Cheng and Liu (\citeyear{ZhangHL2011}) proposed a method for
determining the zero, linear and nonlinear components in partially
linear models. Their method is a regularization method in the
smoothing spline ANOVA framework that is closely related to the COSSO.
They obtained the rate of convergence of their proposed estimator. They
also showed that their method is selection consistent in the special
case of tensor product design. But their approach requires tuning of
four penalty parameters, which may be difficult to implement in
practice. Huang, Wei and Ma (\citeyear{HuangWM2011}) proposed a semiparametric regression
pursuit method for estimating $S_1$ and $S_2$. They embedded partially
linear models into model (\ref{GamA}). By approximating the
nonparametric components using spline series expansions, they
transformed the problem of model specification into a group variable
selection problem. They then used the 2-norm group MCP to determine the
linear and nonlinear components. They showed that, under suitable
conditions, the proposed approach is consistent in estimating the
structure of (\ref{PLM1}), meaning that it can correctly determine
which covariates have a linear effect and which do not with high
probability.

\subsection{Varying Coefficient Models} Consider the linear
varying coefficient model
%
\begin{eqnarray}\label{VcmA}
y_{i}(t_{ij})=\sum_{k=1}^{p}x_{ik}(t_{ij})\beta_{k}(t_{ij})+\epsilon
_{i}(t_{ij}),\nonumber\\
\eqntext{\mbox{$i=1,\ldots,n$, $j=1,\ldots,n_{i}$,}}
\end{eqnarray}
where $y_{i}(t)$ is the response variable for the $i$th subject at time
point $t \in T$ with $T$ being the time interval on which the
measurements are taken, $\epsilon_{i}(t)$ is the error term,
$x_{ik}(t)$'s are time-varying covariates, $\beta_{k}(t)$ is the
corresponding smooth coefficient function. Such a model is useful in
investigating the time-dependent effects of covariates on responses
measured repeatedly. One well-known example is longitudinal data
analysis (Hoover et al., \citeyear{HRWY}) where the response for the $i$th
experimental subject in the study is observed $n_{i}$ occasions, the\vadjust{\goodbreak}
set of observations at times $\{t_{ij}\dvtx j=1,\ldots, n_{i}\}$ tends to
be correlated. Another important example is the functional response
models (Rice, \citeyear{Rice2004}), where the response $y_{i}(t)$ is a smooth real
function, although only $y_{i}(t_{ij})$, $j=1, \ldots, n_{i}$ are
observed in practice. Wang, Chen and Li (\citeyear{WangCL2007}) and Wang and Xia (\citeyear{WangXia})
considered the use of group LASSO and SCAD methods for model selection
and estimation in (\ref{VcmA}). Xue, Qu and Zhu (\citeyear{LAZ2010}) applied the
2-norm SCAD method for variable selection in generalized linear
varying-coefficient models and considered its selection and estimation
properties. These authors obtained their results in the models with
fixed dimensions. Wei, Huang and Li (\citeyear{WHL2010}) studied the properties of
the group LASSO and adaptive group LASSO for (\ref{VcmA}) in the $p
\gg
n$ settings. They showed that, under the sparse Riesz condition and
other regularity conditions, the group LASSO selects a model of the
right order of dimensionality, selects all variables with coefficient
functions whose $\ell_{2}$ norm is greater than a certain threshold
level and is estimation consistent. They also proved that the adaptive
group LASSO can correctly select important variables with high
probability based on an initial consistent estimator.\vspace*{-2pt}

\subsection{Seemingly Unrelated Regressions and Multi-Task
Learning} Consider $T$ linear regression models
\[
\by_t=X_t\bbeta_t + \bveps_t,\quad t=1, \ldots, T,
\]
where $\by_t$ is an $n\times1$ response vector, $X_t$ is an $n\times
p$ design matrix, $\bbeta_t$ is a $p \times1$ vector of regression
coefficients and $\bveps_t$ is an $n\times1$ error vector. Assume that
$\bveps_1, \ldots, \bveps_T$ are independent and identically
distributed with mean zero and covariance matrix $\Sigma$. This model
is referred to as the seemingly unrelated regressions (SUR) model
(Zellner, \citeyear{Zellner1962}). Although each model can be estimated separately based
on least squares method, it is possible to improve on the estimation
efficiency of this approach. Zellner (\citeyear{Zellner1962}) proposed a method for
estimating all the coefficients simultaneously that is more efficient
than the single-equation least squares estimators. This model is also
called a multi-task learning model in machine learning (Caruana, \citeyear{CaruanaR1997};
Argyriou, Evgeniou and Pontil, \citeyear{AEP}).

Several authors have considered the problem of variable selection based
on the criterion
\[
\frac{1}{2T}\sum_{t=1}^T \frac{1}{n}\|\by_t-X_t\bbeta_t\|_2^2 +
\lam
\sum_{j=1}^p \Biggl(\sum_{t=1}^T \beta_{tj}^2\Biggr)^{1/2}.
\]
This is a special case of the general group LASSO criterion. Here the\vadjust{\goodbreak}
groups are formed by the coefficients corresponding to the $j$th
variable across the regressions. The assumption here is that the $j$th
variable plays a similar role across the tasks and should be selected
or dropped at the same time. Several authors have studied the
selection, estimation and prediction properties of the group LASSO in
this model; see, for example, Bach (\citeyear{Bach2008}), Lounici et al. (\citeyear{LPTG2009}),
Lounici et al. (\citeyear{LPTG2011}) and Obozinski, Wainwright and Jordan (\citeyear{OWJ}) under various
conditions on the design matrices and other regularity conditions.

\subsection{Analysis of Genomic Data} Group selection methods
have important applications in the analysis of high throughput genomic
data---for example, to find genes and genetic pathways that affect a
clinical phenotype such as disease status or survival using gene
expression data. Most phenotypes are the result of alterations in a
limited number of pathways, and there is coordination among the genes
in these pathways. The genes in the same pathway or functional group
can be treated as a group. Efficiency may be improved upon by
incorporating pathway information into the analysis, thereby selecting
pathways and genes simultaneously. Another example is integrative
analysis of multiple genomic datasets. In gene profiling studies,
markers identified from analysis of single datasets often suffer from a
lack of reproducibility. Among the many possible causes, the most
important one is perhaps the relatively small sample sizes and hence
lack of power of individual studies. A cost-effective remedy to the
small sample size problem is to pool and analyze data from multiple
studies of the same disease. A generalized seemingly unrelated
regressions model can be used in this context, where a group structure
arises naturally for the multiple measurements for the same gene across
the studies. Some examples of using group selection methods in these
applications include Wei and Li (\citeyear{WeiL2007}), Jacob, Obozinski and Vert (\citeyear{JOV}), Ma and
Huang (\citeyear{MaH2009}), Ma, Huang and Moran (\citeyear{MaHM2009}),
Ma, Huang and Song (\citeyear{MHSong2010}), Ma et al. (\citeyear{MHetal2011}),
Pan, Xie and Shen (\citeyear{PXS}) and Peng et al. (\citeyear{Pengetal}), among others.

\subsection{Genome Wide Association Studies} Genome wide
association studies (GWAS) are an important method for identifying
disease susceptibility genes for common and complex diseases. GWAS
involve scanning hundreds to thousands of samples, often as
case-control samples, utilizing hundreds of thousands of single
nucleotide polymorphism (SNP) markers located\vadjust{\goodbreak} throughout the human
genome. The SNPs from the same gene can be naturally considered as a
group. It is more powerful to select both SNPs and genes simultaneously
than to select them separately. Applications of group selection methods
to genetic association analysis are discussed in Breheny and Huang
(\citeyear{BrehenyHuang1}) and Zhou et al. (\citeyear{ZhouSSL2010}).

\section{Discussion}
In this article, we provide a selective review of several group
selection and bi-level selection methods. While considerable progress
has been made in this area, much work remains to be done on a number of
important issues. Here we highlight some of them that require further
study in order to better apply these methods in practice.

\subsection{Penalty Parameter Selection} In any penalization
approach for variable selection, a difficult question is how to
determine the penalty parameters. This question is even more difficult
in group selection methods. Widely used criterions, including the AIC
(Akaike, \citeyear{Akaike1973}) and BIC (Schwarz, \citeyear{Sch1978}), require the estimation of the
error variance and degrees of freedom. For the group LASSO, Yuan and
Lin (\citeyear{YuanL2006}) proposed an estimate of the degrees of freedom, but it
involves the least squares estimator of the coefficients, which is not
feasible in $p \gg n$ models. The problem of variance estimation has
not been studied systematically in group selection models. It is
possible to use $K$-fold cross validation, which does not require
estimating the variance or the degrees of freedom. However, to our
knowledge, there have been no rigorous analyses of this procedure in
group selection settings. Recently, Meinshausen and B\"{u}hlmann (\citeyear{MeinshausenB10})
proposed stability selection for choosing penalty parameters based on
resampling. This is a general approach
and is applicable to the group selection methods discussed here.
Furthermore, it does not require estimating the variance or
the degrees of freedom. It would be interesting to apply this method to
group selection and compare it with the existing methods in group
selection problems.

\subsection{Theoretical Properties} Currently, most
theoretical results concerning selection, estimation and prediction on
group selection methods in $p \gg n$ settings are derived for the group
LASSO in the context of linear regression. These results provide
important insights into the behavior of the group LASSO. However, they
are obtained for a given range of the penalty parameter. It is not
clear whether, if the penalty parameter is selected using a data-driven
procedure, such as cross validation, these results still hold. It is
clearly of practical interest to confirm the estimation and prediction
properties of group LASSO if the penalty parameter is selected using
such a procedure. For concave selection methods, we considered the
selection property of the global 2-norm group MCP solutions. Although
global results shed some light on the properties of these methods, it
is more relevant to investigate the properties of the local solutions,
such as those obtained based on the group coordinate descent algorithm.
Therefore, much work is needed to understand the theoretical properties
of various concave group selection methods and compare their
performance with the group LASSO.

\subsection{Overlapping Groups} In this article, we only
considered the case where there is no overlapping among the groups.
However, in many applied problems, overlapped groups arise naturally.
For example, in genomic data analysis involving genes and pathways,
many important genes belong to multiple pathways. Jacob, Obozinski and
Vert (\citeyear{JOV}) proposed an extended group LASSO method for selection with
overlapping groups. With their method, it is possible to select one
variable without selecting all the groups containing it. Percival
(\citeyear{Pervival2011}) studied the theoretical properties of the method of Jacob, Obozinski and
Vert (\citeyear{JOV}). Liu and Ye (\citeyear{LY2010}) proposed an algorithm for solving the
overlapping group LASSO problem. Zhao, Rocha and Yu (\citeyear{ZhaoRY2009}) considered the
problem of overlapping groups in the context of composite absolute
penalties. The results of Huang et al. (\citeyear{HuangMXZ}) on the selection
consistency of the 1-norm group bridge allow overlapping among groups
under the assumption that the extent of overlapping is not large.
However, in general, especially for concave group selection methods,
this question has not been addressed.

\appendix
\section*{Appendix}\label{appendix}

\begin{lemma}
\label{LemChi2} Let $\chi_k^2$ be a random variable with chi-square
distribution with $k$ degrees of freedom. For $t > 1$,
$\rP(\chi^2_k\ge k t) \le h(t,k),$ where $h(t,k)$ is defined in
(\ref{Defh}).
\end{lemma}
This lemma is a restatement of the exponential inequality for
chi-square distributions of Laurent and Massart (\citeyear{LaurentM2010}).

\begin{pf*}{Proof of Theorem \ref{ThmA}}
Since $\hbbeta^o$ is
the oracle least squares estimator, we have $\hbbeta_j^o=0$ for $j
\notin S$ and
%
\begin{equation}\label{KKTa}
-X_j'(\by-X\hbbeta^o)/n =0\quad\forall\! j \in S.
\end{equation}
If $\|\hbbeta_j^o\|_2/\sqrt{d_j} \ge\gam\lam$, then by the definition
of the MCP, $\rho'(\|\hbbeta_j^o\|_2; \sqrt{d_j}\lam, \gam)=0$. Since
$c_{\min} > 1/\gam$, the criterion (\ref{gM1}) is strictly convex. By
the KKT conditions, the equality $\hbbeta(\lam, \gam)=\hbbeta^o$ holds
in the intersection of the events
%
\begin{eqnarray}\label{KKTb}
\Omega_1(\lam) &=& \Bigl\{ \max_{j \notin S}
\|n^{-1}X_j'(\by-X\hbbeta^o) \|_2/\sqrt{d_j} \nonumber\\
&&\hspace*{127.4pt}\le\lam\Bigr\}
\quad\mbox{and}\\
 \Omega_2(\lam) &=& \Bigl\{ \min_{j \in S}
\|\hbbeta_j^o\|_2/\sqrt{d_j} \ge\gam\lam\Bigr\}.\nonumber
\end{eqnarray}

We first bound $1-\rP(\Omega_1(\lam))$.
Let $\hbbeta_{S}=(\hbbeta_j, j \in S)'$. 
By (\ref{KKTa})
and using $\by=X_{S}\bbeta_{S}^o+\bveps$,
%
\begin{equation}\label{KKTc}
\hbbeta_{S}^o = \Sigma_{S}^{-1}X_{S}'\by/n= \bbeta_{S}^o
+\Sigma_{S}^{-1}X_{S}'\bveps/n.
\end{equation}
It follows that \( n^{-1}X_j'(\by-X\hbbeta^o)=n^{-1}
X_j'(I_n-P_{S})\bveps, \) where
$P_{S}=n^{-1}X_{S}\Sigma_{S}^{-1}X_{S}'$. Because
$X_j'X_j=I_{d_j}$,\break
$\|X_j'(I_n-P_{S})\bveps\|_2^2/\sigma^2$ is distributed as a $\chi^2$
distribution with $d_j$ degrees of freedom. We have, for
$n\lam^2/\sigma^2 \ge1$,
%
\begin{eqnarray}\label{Aeq1}
&&1-\rP(\Omega_1(\lam))\nonumber\\
 &&\quad= \rP\Bigl(\max_{j \notin S}\|
n^{-1/2}X_j'(I_n-P_{S})\bveps\|_2^2/(d_j\sigma^2)\nonumber\\
&&\hspace*{138.5pt}\quad > n \lam^2 /\sigma
^2 \Bigr)\nonumber\\
&&\quad\le \sum_{j \notin S} \rP\bigl(\|n^{-1/2}X_j'(I_n-P_{S})\bveps
\|_2^2/\sigma^2 \nonumber\\[-8pt]\\ [-8pt]
&&\hspace*{113.5pt}>d_j n \lam^2 /\sigma^2 \bigr)\nonumber\\
&&\quad\le \sum_{j \notin S} h(n\lam^2/\sigma^2, d_j) \nonumber\\
&&\quad\le (J-|S|) h \bigl(n\lam^2/\sigma^2, d_{\min}(S^c) \bigr)\nonumber\\
&&\quad= \eta_{1n}(\lam),\nonumber
\end{eqnarray}
where we used Lemma \ref{LemChi2} in the third line.

Now consider $\Omega_2$. Recall $\beta_*^o=\min_{j \in
S}\|\bbeta_j^o\|_2/\sqrt{d_j}$. If
$\|\hbbeta_j^o-\bbeta_j^o\|_2/\sqrt{d_j}\le\beta_*^o-\gam\lam$
for all
$j\in S$, then $\min_{j \in S}\|\hbbeta_j^o\|_2/\sqrt{d_j} \ge
\gam\lam$. This implies
\[
\label{Aeq2a} 1-\rP(\Omega_2(\lam)) \le\rP\Bigl(\max_{j \in
S}\|\hbbeta_j^o-\bbeta_j^o\|_2/\sqrt{d_j}
> \beta_*^o-\gam\lam\Bigr).
\]
Let $A_j$ be a $d_j \times d_S$ matrix with a $d_j \times d_j$ identity
matrix $I_{d_j}$ in the $j$th block and $0$'s elsewhere. Then \(
n^{1/2}(\hbbeta_j^o-\bbeta_j^o)=
n^{-1/2}A_j\Sigma_{S}^{-1}X_{S}'\bveps. \) Note that
\begin{eqnarray*}
&&\|n^{-1/2}A_j\Sigma_{S}^{-1}X_{S}'\bveps\|_2\\
 &&\quad\le\|A_j\|_2
\|\Sigma_S^{-1/2}\|_2 \|n^{-1/2}\Sigma_S^{-1/2}X_S'\bveps\|_2\\
 &&\quad\le
c_1^{-1/2}\|n^{-1/2}\Sigma_S^{-1/2}X_S'\bveps\|_2
\end{eqnarray*}
and $\|n^{-1/2}\Sigma_S^{-1/2}X_S'\bveps\|_2^2/\sigma^2$ is distributed
as a $\chi^2$ distribution with $|S|$ degrees of freedom. Therefore,
similar to (\ref{Aeq1}), we have, for $c_1 n
(\beta_*^o-\gam\lam)^2/\sigma^2 > 1$,
%
\begin{eqnarray}\label{Aeq2}
&&1-\rP(\Omega_2(\lam))\nonumber\\
 &&\quad=\rP\Bigl(\max_{j \in S} n^{-1/2}
\|A_j\Sigma_{S}^{-1}X_{S}'\veps\|_2/\sqrt{d_j}\nonumber\\
&&\hspace*{100pt}> \sqrt{n}(\beta_*^o-\gam\lam) \Bigr) \nonumber\\
&&\quad\le \rP\Bigl(\max_{j \in
S}\|n^{-1/2}\Sigma_S^{-1/2}X_S'\bveps\|_2^2/(d_j\sigma^2)\\
&&\hspace*{84pt}> c_1 n (\beta_*^o-\gam\lam)^2/\sigma^2 \Bigr) \nonumber\\
&&\quad\le |S| h \bigl(c_1 n (\beta_*^o-\gam\lam)^2/\sigma^{2}, d_{\min}(S) \bigr)\nonumber\\
&&\quad= \eta_{2n}(\lam).\nonumber
\end{eqnarray}
Combining (\ref{Aeq1}) and (\ref{Aeq2}), we have
\begin{eqnarray*}
\rP\bigl(\hbbeta(\lam, \gam) \neq\hbbeta^o\bigr) &\le& 1-\rP(\Omega_1(\lam
)) + 1-
\rP(\Omega_2(\lam))\\
 &\le&\eta_{1n}(\lam) + \eta_{2n}(\lam).
\end{eqnarray*}
This completes the proof.
\end{pf*}

For any $B \subset\{1, \ldots, J\}$ and $m \ge1$, define
%
\begin{eqnarray}\label{zetanorm}
&&\hspace*{20pt}\zeta(\upsilon; m,B)\nonumber\\
&&\hspace*{20pt}\quad=\max\biggl\{
\frac{\|(P_{A}-P_{B})\upsilon\|_{2}}{(mn)^{1/2}}\dvtx \\
&&\hspace*{20pt}\quad \hspace*{36pt}B \subseteq A
\subseteq\{1, \ldots, J\}, d_{A}=m+d_{B} \biggr\}\nonumber
\end{eqnarray}
for $\upsilon\in\mathbb{R}^{n}$, where $P_{A}=X_A(X_A'X_A)^{-1}X_A'$
is the orthogonal projection from $\mathbb{R}^{n}$ to the span of
$X_{A}$.

\begin{lemma}
\label{LemB} Suppose $\xi n\lam^2> \sigma^2 d_{\max}$. We have
\begin{eqnarray*}
&&\rP\bigl(2\sqrt{c^*d_s}\zeta(\by;m, S)> \lam\bigr)\\
&&\quad \le(J-|S|)^m \frac{e^m}{m^m}
\exp(-m\xi n\lam^2/16).
\end{eqnarray*}
\end{lemma}

\begin{pf}
For any $A \supseteq S$, we have
$(P_A-P_S)\cdot X_{S}\bbeta_{S}=0$. Thus
\[
(P_A-P_S)\by=(P_A-P_{A_{S}})(X_{S}\bbeta_{S}+\veps)= (P_A-P_S)\veps.
\]
Therefore,
\begin{eqnarray*}
&&\rP\bigl(2\sqrt{c^*d_s }\zeta(\by;m, S)> \lam\bigr)\\
 &&\quad=\rP\Bigl(\max_{A
\supseteq S,
|A|-|S|=m} \|(P_A-P_{S})\veps\|^2/\sigma^2> \xi m n\lam^2\Bigr).
\end{eqnarray*}
Since $P_A-P_{S}$ is a projection matrix,
$\|(P_A-P_{S})\veps\|^2/\allowbreak \sigma^2 \sim\chi^2_{m_A}$, where $m_A =
\sum_{j \in A-S, A \supseteq S} d_j \le m d_{\max}$.\break Since there are
${J-|S| \choose m}$ ways to choose $A$ from $\{1, \ldots,\break J\}$, we have
\begin{eqnarray*}
&&\rP\bigl(2\sqrt{c^*d_s}\zeta(\by;m, S)> \lam\bigr) \\
&&\quad\le{J-|S| \choose m}
\rP(\chi_{md_{\max}}^2 > \xi mn \lam^2).
\end{eqnarray*}
This and Lemma \ref{LemChi2} imply that
\begin{eqnarray*}
&&\rP\bigl(2\sqrt{c^*d_s}\zeta(\by;m, S)> \lam\bigr)\\
&&\quad\le {J-|S| \choose m} h(\xi n\lam^2/d_{\max}, m d_{\max})\\
&&\quad\le (J-|S|)^m \frac{e^m}{m^m} h(\xi n\lam^2/d_{\max}, m d_{\max
}). 
\end{eqnarray*}
Here we used the inequality ${J-|S| \choose m} \le e^m(J-|S|)^m/\break m^m$.
This completes the proof.
\end{pf}

Define $T$ as any set that satisfies
\begin{eqnarray*}
&&S\cup\{j\dvtx \|\hat{\bbeta}_{j}\|_{2} \ne0\} \\
&&\quad\subseteq T \subseteq
S\cup
\bigl\{j\dvtx
n^{-1}X_{j}^{\prime}(\by-X\hat{\bbeta})\\
&&\quad\phantom{\subseteq T \subseteq
S\cup
\bigl\{j\dvtx}=\dot{\rho}\bigl(\|\hat{\bbeta
}_{j}\|_{2};
\sqrt{d_{j}}\lambda, \gamma\bigr)\\
&&\hspace*{104pt}{}\cdot{\sqrt{d_{j}}\hat{\bbeta}_{j}}/{\|\hat{\bbeta}_{j}\|_{2}}\bigr\}.
\end{eqnarray*}

\begin{lemma}
\label{LemC} Suppose that $X$ satisfies the\break $\operatorname{SRC}(d^*, c_*, c^*)$,
$d^*\ge(K_*+1)|S|d_s $, and
$\gam\ge c_*^{-1}\cdot\sqrt{4+c_*/c^*}$.
Let $m_*=K_*|S|$. 
Then for any $\by\in\real^{n}$
with $\lam\ge2\sqrt{c^*d_s} \zeta(\by; m_*, S)$,
we have
\[
|T| \le(K_{*}+1)|S|.
\]
\end{lemma}

\begin{pf}
This lemma can be proved along the line of
the proof of Lemma 1 of Zhang (\citeyear{Zhang2010}) and is omitted.
\end{pf}

\begin{pf*}{Proof of Theorem \ref{ThmB}}
By Lemma \ref{LemC}, in
the event
%
\begin{equation}\label{zetaConA}
2\sqrt{c^*d_{\max}(S)}  \zeta(\by;m_*, S)\le
\lam,
\end{equation}
we have $|T| \le(K_*+1)|S|$. Thus in event (\ref{zetaConA}), the
original model with $J$ groups reduces to a model with at most
$(K_*+1)|S|$ groups. In this reduced model, the conditions of Theorem
\ref{ThmB} imply that the~conditions\vadjust{\goodbreak} of Theorem \ref{ThmA} are
satisfied. By\break Lemma~\ref{LemB},
%
\begin{equation}\label{PboundAA}
\hspace*{20pt}\rP\bigl(2\sqrt{c^*d_{\max}(S)} \zeta(\by;m_*, S)>
\lam\bigr)
\le\eta_{3n}(\lam).
\end{equation}
Therefore, combining (\ref{PboundAA}) and Theorem \ref{ThmA}, we have
\[
\label{PboundB} \rP\bigl( \hbbeta(\lam, \gam) \neq\hbbeta^o\bigr) \le
\eta_{1n}(\lam) + \eta_{2n}(\lam)+\eta_{3n}(\lam).
\]
This proves Theorem \ref{ThmB}.
\end{pf*}

\section*{Acknowledgments}
We wish to thank two anonymous reviewers, the Associate Editor and
Editor for their helpful comments. In particular, we are extremely
grateful to one reviewer for providing valuable and detailed comments
and for pointing out the work of Laurent and Massart (\citeyear{LaurentM2010}) to us,
which have led to substantial improvements in the paper. The research
of Huang is partially supported by NIH Grants R01CA120988, R01CA142774
and NSF Grant DMS-08-05670. The research of Ma is partially supported by
NIH Grants R01CA120988 and R01CA142774.


%
\end{document}